\documentclass[12pt]{article}
\usepackage{amsmath,amssymb,amscd,verbatim,amsthm}
\usepackage{latexsym, epsfig, graphpap, graphics, mathrsfs}
\usepackage[varg]{pxfonts}
\usepackage{mathrsfs}

\begin{document}
 \noindent{\bf{On strongly almost lacunary statistical $A$-convergence defined by Musielak-Orlicz function}}

\vskip 0.5 cm

\noindent Ekrem Sava$\c{s}$

\noindent Istanbul Commerce University, 34840 Istanbul, Turkey

\noindent E-mail : ekremsavas@yahoo.com
\vskip 0.5 cm

\noindent Stuti Borgohain $^*${\footnote {The work of the
authors was carried under the Post Doctoral Fellow under National Board of Higher 
Mathematics, DAE, project No. NBHM/PDF.50/2011/64}}

\noindent Department of Mathematics

\noindent Indian Institute of Technology, Bombay

\noindent Powai:400076, Mumbai, Maharashtra; INDIA.

\noindent E-mail:stutiborgohain$@$yahoo.com
\vskip 1 cm

\noindent{\footnotesize {\bf{Abstract:}} We study some new strongly almost lacunary statistical $A$-convergent sequence space of order $\alpha$ defined by a Musielak-Orlicz function. We also give some inclusion relations between the newly introduced class of sequences with the spaces of strongly almost lacunary $A$-convergent sequence of order $\alpha$. Moreover we have examined some results on Musielak-Orlicz function with respect to these spaces.

 
\vskip 0.3 cm

\noindent{\bf{AMS Classification No:} 40A05; 40A25; 40A30;
40C05.}}

\vskip 1 cm

\section{Introduction}

The concept of statistical convergence was initially introduced by Fast \cite{Fast}, which is closely related to the concept of natural density or asymtotic density of subsets of the set of natural numbers $N$. Later on, it was studied as asummability method by Fridy \cite{Fridy}, Fridy and Orhan \cite{Orhan}, Freedman and Sember \cite{Freedman}, Schoenberg \cite{Schoenberg}, Malafosse and Rako$\check{c}$evi$\acute{c}$ \cite{Malafossa} and many more mathematicians. Moreover, in recent years, generalizations of statistical convergence have appeared in the study of strong integral summability and the structure of ideals of bounded continuous functions on locally compact spaces. Also, statistical convergence is closely related to the concept of convergence in probabilty.  

By the concept of almost convergence, we have a sequence $x=(x_k) \in \ell_\infty$ if its Banach limit coincides. The set $\hat c$ denotes set of all almost convergent sequences .  Lorentz \cite{Lorentz} proved that,

$$\hat c =\{ x \in \ell_\infty: \displaystyle\lim_m t_{mn}(x) \mbox{~exist uniformly in~} n \},$$

where 

$$t_{mn}(x)= \frac{x_n+x_{n+1}+....+x_{n+m}}{m+1}.$$

Similarly, the space of strongly almost convergent sequence was defined as, $[\hat c]=\{x \in \ell_\infty: \displaystyle\lim_m t_{m,n}(\vert x -Le \vert) \mbox{~exists uniformly in ~} n \mbox{~ for some ~} L \},$ where, $e=(1,1,...)$.  (see Maddox \cite{Maddox})\\

A lacunary sequence is defined as an increasing integer sequence $\theta = (k_r)$ such that $k_0=0$ and $h_r=k_r-k_{r-1} \rightarrow \infty$ as $r \rightarrow \infty.$\\

\noindent Note: Throughout this paper, the intervals determined by $\theta$ will be denoted by $J_r=(k_{r-1}, k_r]$ and the ratio $\frac{k_r}{k_{r-1}}$ will be defined by $\phi_r$. \\

\section{Preliminary concepts}

Let $0<\alpha \leq 1$ be given. The sequence $(x_k)$ is said to be statistically convergent of order $\alpha$ if there is a real number $L$ such that,

$$\displaystyle\lim_{n \rightarrow \infty} \frac{1}{n^\alpha} \vert \{ k \leq n : \vert x_k -L \vert \geq \varepsilon \} \vert =0,$$

for every $\varepsilon >0$. In this case, we write $S^\alpha-\displaystyle\lim x_k=L$. The set of all statistically convergent sequences of order $\alpha$ will be denoted by $S^\alpha$.\\

For any lacunary sequence $\theta= (k_r)$, the space $N_\theta$ defined as, (Freedman et al.\cite{Freedman})

$$N_\theta=\left\{(x_k): \displaystyle\lim_{r \rightarrow \infty} h_r^{-1} \displaystyle\sum_{k \in J_r} \vert x_k - L \vert =0, \mbox{~for some~} L \right\}.$$

The space $N_\theta$ is a $BK$ space with the norm,

$$\Vert (x_k) \Vert_\theta= \displaystyle\sup_r h_r^{-1} \displaystyle\sum_{k \in J_r} \vert x_k \vert.$$

Let $\theta=(k_r)$ be a lacunary sequence and $0<\alpha \leq 1$ be given. The sequence $x=(x_k) \in w$ is said to be $S_\theta^\alpha$-statistically convergent (or lacunary statistically convergent sequence of order $\alpha$) if there is a real number $L$ such that

$$\displaystyle\lim_{r \rightarrow \infty} \frac{1}{h_r^\alpha} \vert \{ k \in I_r : \vert x_k -L \vert \geq \varepsilon \} \vert =0,$$

where $I_r =(k_{r-1}, k_r]$ and $h_r^\alpha$ denotes the $\alpha$-th power $(h_r)^\alpha$ of $h_r$, that is, $h^\alpha=(h_r^\alpha)=(h_1^\alpha,h_2^\alpha,...h_r^\alpha,...)$. We write $S_\theta^\alpha-\displaystyle\lim x_k =L$. The set of all $S_\theta^\alpha$-statistically convergent sequences will be denoted by $S_\theta^\alpha$. \\

By an Orlicz function , we mean a function $M: [0,\infty )\rightarrow [0,\infty )$, which is continuous, non-decreasing and convex with $M(0) = 0, M(x)>0$, for $x>0$ and $M(x)\rightarrow \infty$, as $x\rightarrow \infty$.\\

The idea of Orlicz function is used to construct the sequence space, (see Lindenstrauss and Tzafriri \cite{Lindenstrauss}),

$$\ell_M=\left\{ (x_k) \in w: \displaystyle\sum_{k=1}^\infty M \left(\frac{\vert x \vert}{\rho}\right) < \infty, \mbox{~for some~} \rho>0 \right\}. $$

This space $\ell_M$ with the norm,

$$ \Vert x \Vert = \mbox{inf}\left\{ \rho>0: \displaystyle\sum_{k=1}^\infty M \left(\frac{\vert x_k \vert}{\rho}\right) \leq 1 \right\}$$

becomes a Banach space which is called an Orlicz sequnce space.\\

Musielak \cite{Musielak} defined the concept of Musielak-Orlicz function as $\mathscr{M}=(M_k)$. A sequence $\mathscr{N}=(N_k)$ defined by 

$$N_k(v)=\sup \{ \vert v \vert u -M_k(u): u \geq 0 \}, k=1,2,..$$

is called the complementary function of a Musielak-Orlicz function $\mathscr{M}$. The Musielak-Orlicz sequence space $t_\mathscr{M}$ and its subspace $h_\mathscr{M}$ are defined as follows:

$$t_\mathscr{M}=\{ x \in w: I_\mathscr{M}(cx) < \infty \mbox{~for some ~} c>0\},$$
$$h_\mathscr{M}=\{x \in w:I_\mathscr{M}(cx) < \infty, \forall c >0\},$$

where $I_\mathscr{M}$ is a convex modular defined by,

$$I_\mathscr{M}(x) = \displaystyle\sum_{k=1}^\infty M_k(x_k), x=(x_k) \in t_\mathscr{M}.$$

It is considered $t_\mathscr{M}$ equipped with the Luxemberg norm

$$\Vert x \Vert = \inf \left\{ k>0: I_\mathscr{M}\left(\frac{x}{k}\right) \leq 1 \right\}$$

or equiped with the Orlicz norm

$$\Vert x \Vert^0 = \inf \left\{ \frac{1}{k} (1+I_\mathscr{M}(kx)): k>0 \right\}.$$

A Musielak-Orlicz function $(M_k)$ is said to satisfy $\Delta_2$-condition if there exist constants $a, K>0$ and a sequence $c=(c_k)_{k=1}^\infty \in \ell_+^1$ ( the positive cone of $\ell^1$) such that the inequality

$$M_k(2u) \leq KM_k(u)+c_k$$

holds for all $k \in N$ and $u \in R_+$, whenever $M_k(u) \leq a$.

If $A=(a_{nk})_{n,k=1}^\infty$ is an infinite matrix, then $Ax$ is the sequence whose $n$th term is given by $A_n(x)=\displaystyle\sum_{k=1}^\infty a_{nk}x_k$.\\

We consider a sequence $x=(x_k)$ which is said to be strongly almost lacunary statistical $A$-convergent of order $\alpha$ (or $S_\theta^\alpha(A, \mathscr{M}, (s))$-statistically convergent) if ,

$$\displaystyle\lim_{r \rightarrow  \infty} \frac{1}{h_r^\alpha} \left\vert \left\{ k \in I_r : \displaystyle\sum_{k \in I_r} \left( M_k \left(\frac{\vert t_{km} (A_k(x)-L) \vert}{\rho^{(k)}} \right)\right)^{(s_k)} \geq \varepsilon \right\}\right \vert=0, \mbox{~uniformly in ~} m, $$

where $I_r = (k_{r-1}, k_r]$ and $h_r^\alpha$ denotes the $\alpha$-th power $(h_r^\alpha)$ of $h_r$, that is, $h^\alpha=(h_r^\alpha)=(h_1^\alpha, h_2^\alpha,...h_r^\alpha,...)$ and $\mathscr{M}=(M_k)$ is a Musielak-Orlicz function.\\

Also we have introduced the space of strongly almost lacunary $A$-convergent sequences with respect to Musielak-Orlicz function 
$\mathscr{M}=(M_k)$ as follows:

$$\hat N_\theta^\alpha(A, \mathscr{M}, (s))=\left\{ (x_k) :\displaystyle\lim_{r \rightarrow \infty} \frac{1}{h_r^\alpha} \displaystyle\sum_{k \in I_r} \left(M_k \left( \frac{\vert t_{km} (A_k(x)-L) \vert }{\rho^{(k)}} \right)\right) ^{(s_k)} =0, \mbox{~for some ~} L \mbox{~and~} \rho^{(k)} >0 \right\}.$$

We give some inclusion relations between the sets of $S_\theta^\alpha(A, \mathscr{M}, (s))$-statistically convergent sequences and strongly almost lacunary $A$-convergent sequence space $\hat N_\theta^\alpha(A, \mathscr{M}, (s))$. Also some results defined by Musielak-Orlicz function are studied with respect to these sequence spaces.
 
\section{Main Results}

\textbf{Theorem 3.1.} Let $\alpha, \beta \in (0,1]$ be real numbers such that $\alpha\leq \beta$, $\mathscr{M}$ be a Musielak-Orlicz function and $\theta=(k_r)$ be a lacunary sequence, then $\hat N_\theta^\alpha(A,\mathscr{M},(s)) \subset \hat S_\theta^\beta$.  \\

Proof: Let $x \in \hat N_\theta^\alpha(A, \mathscr{M}, (s))$.\\

For $\varepsilon >0$  given, let us denote $\sum_1$ as the sum over $k \in I_r$, $\vert t_{km}(A_k(x)-L) \vert \geq \varepsilon.$ and $\sum_2$ denote the sum over $k \in I_r$, $\vert t_{km}(A_k(x)-L) \vert < \varepsilon$ respectively.\\ 

As $h_r^\alpha \leq h_r^\beta$ for each $r$, we may write,\\

$\frac{1}{h_r^\alpha} \displaystyle\sum_{k \in I_r} \left[M_k \left(\frac{\vert t_{km}(A_k(x)-L) \vert}{\rho^{(k)}} \right)\right]^{s_k}$\\

$=\frac{1}{h_r^\alpha} \left[\sum_1 \left[M_k\left(\frac{\vert t_{km}(A_k(x)-L) \vert}{\rho^{(k)}} \right)\right]^{s_k} +\sum_2 \left[M_k\left(\frac{\vert t_{km}(A_k(x)-L) \vert}{\rho^{(k)}} \right)\right]^{s_k}\right]$\\

$\geq \frac{1}{h_r^\beta} \left[\sum_1 \left[M_k\left(\frac{\vert t_{km}(A_k(x)-L) \vert}{\rho^{(k)}} \right)\right]^{s_k} +\sum_2 \left[M_k\left(\frac{\vert t_{km}(A_k(x)-L) \vert}{\rho^{(k)}} \right)\right]^{s_k}\right]$\\

$\geq \frac{1}{h_r^\beta} \sum_1\left[ M_k\left(\frac{\varepsilon}{\rho^{(k)}}\right)\right]^{s_k}$\\

$\geq \frac{1}{h_r^\beta} \sum_1 \mbox{min}([M_k(\varepsilon_1)]^h, [M_k(\varepsilon_1)]^H), \varepsilon_1=\frac{\varepsilon}{\rho^{(k)}}$\\

$\geq \frac{1}{h_r^\beta} \vert \{k \in I_r : \vert t_{km}(A(x)) -L \vert \geq \varepsilon \} \vert \mbox{min}([M_k(\varepsilon_1)]^h, [M_k(\varepsilon_1)]^H.$ 

As $x \in \hat N_\theta^\alpha(A, \mathscr{M}, (s))$, the left hand side of the above inequality tends to zero as $ r \rightarrow \infty$.\\

Therefore, the right hand side of the above inequality tends to zero as $r \rightarrow \infty$, hence $x \in \hat S_\theta^\beta$.\\

\textbf{Corollary 3.2.}  Let $0 < \alpha \leq 1$, $\mathscr{M}$ be a Musielak-Orlicz function and $\theta=(k_r)$ be a lacunary sequence, then

$$\hat N_\theta^\alpha(A, \mathscr{M}, (s)) \subset \hat S_\theta^\alpha.$$

\textbf{Theorem 3.3.}  Let $\mathscr{M}$ be a Musielak-Orlicz function, $x=(x_k)$ be a bounded sequence and $\theta=(k_r)$ be a lacunary sequence. If $\displaystyle\lim_{r \rightarrow \infty} \frac{h_r}{h_r^\alpha}=1$,then $x \in \hat S_\theta^\alpha \Rightarrow x \in \hat N _\theta^\alpha(A, \mathscr{M}, (s))$.\\

Proof: Suppose that $x=(x_k)$ be a bounded sequence that is $x \in \ell_\infty$ and $\hat S_\theta^\alpha-\mbox{lim} x_k =L$.\\

As $x \in \ell_\infty$, then there is a constant $T>0$ such that $\vert x_k \vert \leq T$. Given $\varepsilon > 0$, we have,\\

$\frac{1}{h_r^\alpha} \displaystyle\sum_{k \in I_r} \left[M_k\left(\frac{\vert t_{km}(A_k(x)-L) \vert}{\rho^{(k)}}\right)\right]^{s_k}$\\

$\frac{1}{h_r^\alpha} \sum_1 \left[M_k\left(\frac{\vert t_{km}(A_k(x)-L) \vert}{\rho^{(k)}}\right)\right]^{s_k} +
\frac{1}{h_r^\alpha} \sum_2 \left[M_k\left(\frac{\vert t_{km}(A_k(x)-L) \vert}{\rho^{(k)}}\right)\right]^{s_k}$\\

$\leq \frac{1}{h_r^\alpha} \sum_1 \mbox{max} \left\{\left[M_k\left(\frac{T}{\rho^{(k)}}\right)\right]^h, \left[M_k\left(\frac{T}{\rho^{(k)}}\right)\right]^H \right\} + \frac{1}{h_r^\alpha} \sum_2 \left[M_k\left(\frac{\varepsilon}{\rho^{(k)}} \right)\right]^{s_k}$\\

$\leq \mbox{max}\{[M_k(K)]^h,[M_k(K)]^H \} \frac{1}{h_r^\alpha} \vert \{ k \in I_r: \vert t_{km}(A_k(x)-L) \vert \geq \varepsilon \} \vert + \frac{h_r}{h_r^\alpha} \mbox{max} \{ [ M_k(\varepsilon_1) ]^h, [M_k(\varepsilon_1)]^H \}, ~~\frac{T}{\rho^{(k)}}=K, \frac{\varepsilon}{\rho^{(k)}}=\varepsilon_1.$\\

Hence, $x \in \hat N_\theta^\alpha(A, \mathscr{M},(s))$. \\

\textbf{Theorem 3.4. } If $\displaystyle\lim s_k >0$ and $x=(x_k)$ is strongly $\hat N_\theta^\alpha(A, \mathscr{M},(s))$-summable to $L$ with respect to the Musielak-Orlicz function $\mathscr{M}$, then $\hat N_\theta^\alpha (A, \mathscr{M}, (s))-\mbox{lim} x_k$ is unique.\\

Proof: Let $\displaystyle\lim s_k =s >0$. Suppose that $\hat N_\theta^\alpha(A, \mathscr{M},(s))-\mbox{lim} x_k =L$, and $\hat N_\theta^\alpha(A, \mathscr{M},(s))-\mbox{lim} x_k =L_1.$ Then,

$$\displaystyle\lim_{r \rightarrow \infty} \frac{1}{h_r^\alpha} \displaystyle\sum_{k \in I_r} \left[ M_k \left( \frac{\vert t_{km}(A_k(x)-L) \vert}{\rho^{(k)}_1} \right) \right] ^{s_k} = 0, \mbox{~for some~} \rho^{(k)}_1>0$$

and 

$$\displaystyle\lim_{r \rightarrow \infty} \frac{1}{h_r^\alpha} \displaystyle\sum_{k \in I_r} \left[ M_k \left( \frac{\vert t_{km}(A_k(x)-L) \vert}{\rho^{(k)}_2} \right) \right] ^{s_k} = 0, \mbox{~for some~} \rho^{(k)}_2>0.$$

Define $\rho^{(k)}=\mbox{max}(2 \rho^{(k)}_1, 2\rho^{(k)}_2)$. As $\mathscr{M}$ is nondecreasing and convex, we have,\\

$\frac{1}{h_r^\alpha} \displaystyle\sum_{k\in I_r} \left[M_k \left( \frac{\vert L - L_1 \vert }{\rho^{(k)}} \right) \right]^{s_k}$\\

$\leq \frac{D}{h_r^\alpha} \displaystyle\sum_{k \in I_r} \frac{1}{2^{s_k}} \left ( \left[ M_k \left( \frac{\vert t_{km}(A_k(x)-L) \vert}{\rho^{(k)}_1} \right)\right] ^{s_k} + \left[ M_k \left( \frac{\vert t_{km}(A_k(x)-L) \vert}{\rho^{(k)}_2} \right)\right] ^{s_k} \right)$\\

$\leq \frac{D}{h_r^\alpha} \displaystyle\sum_{k \in I_r} \left[ M_k \left( \frac{\vert t_{km}(A_k(x)-L) \vert}{\rho^{(k)}_1} \right)\right] ^{s_k} + \frac{D}{h_r^\alpha}\displaystyle\sum_{k \in I_r}\left[ M_k \left( \frac{\vert t_{km}(A_k(x)-L) \vert}{\rho^{(k)}_2} \right)\right] ^{s_k}  \rightarrow 0, (r \rightarrow \infty),$\\

where $\displaystyle\sup_k s_k = H$ and $D=\mbox{max}(1, 2^{H-1})$. Hence,

$$\displaystyle\lim_{r \rightarrow \infty} \frac{1}{h_r^\alpha} \displaystyle\sum_{k \in I_r} \left[M_k\left(\frac{\vert L- L_1 \vert}{\rho^{(k)}} \right)\right]^{s_k}=0.$$

As $\displaystyle\lim_{k \rightarrow \infty} s_k =s$, we have,

$$\displaystyle\lim_{k \rightarrow \infty} \left[M_k \left(\frac{\vert L- L_1 \vert}{\rho^{(k)}} \right)\right]^{s_k} = \left[M_k \left(\frac{ \vert L-L_1\vert}{\rho^{(k)}} \right) \right]^s$$

and so $L = L_1$.Thus the lmit is unique. \\

{\bf{Theorem 3.5. }} Let $A=(a_{mk})$ be an infinite matrix of complex numbers and let $\mathscr{M}=(M_k)$ be a Musielak-Orlicz function satisfying$\Delta_2$-condition . If $x$ is strongly almost lacunary $A$-convergent sequences with respect to $\mathscr{M}$, then $\hat N_\theta^\alpha(A) \subset \hat N_\theta^\alpha(A, \mathscr{M})$.\\

Proof: Let $x \in \hat N_\theta^\alpha(A)$.

Then, $\displaystyle\lim_{r \rightarrow \infty} \frac{1}{h_r^\alpha} \displaystyle\sum_{k \in I_r} \vert t_{km}(A(x)-L) \vert =0$, uniformly in $m$.\\

Let us define two sequences $y$ and $z$ such that,

\[ (\vert t_{km}(A_k(y)-L) \vert) = \left\{ \begin{array}{ll}
(\vert t_{km}(A_k(x)-L) \vert) & \mbox{if $(\vert t_{km}(A_k(x)-L) \vert) >1 $};\\
\theta & \mbox{if $(\vert t_{km}(A_k(x)-L) \vert) \leq 1$}.\end{array} \right. \]

\[ (\vert t_{km}(A_k(z)-L) \vert) = \left\{ \begin{array}{ll}
\theta  & \mbox{if $(\vert t_{km}(A_k(x)-L) \vert) >1 $};\\
(\vert t_{km}(A_k(x)-L) \vert) & \mbox{if $(\vert t_{km}(A_k(x)-L) \vert) \leq 1$}.\end{array} \right. \]\\

Hence, $(\vert t_{km}(A_k(x)-L) \vert)=(\vert t_{km}(A_k(y)-L) \vert) + (\vert t_{km}(A_k(z)-L) \vert).$\\

Also, $(\vert t_{km}(A_k(y )-L) \vert) \leq (\vert t_{km}(A_k(x)-L) \vert)$ and $(\vert t_{km}(A_k(z)-L) \vert) \leq (\vert t_{km}(A_k(x)-L) \vert)$.\\

Since, $\hat N_\theta^\alpha(A)$ is normal, so we have $y,z \in \hat N_\theta^\alpha(A)$.\\

Let $\displaystyle\sup_k M_k (2)=T$\\

Then, \\

$~~\frac{1}{h_r^\alpha} \displaystyle\sum_{k \in I_r} \left[ M_k \left( \frac{ \vert t_{km}(A_k(x)-L) \vert }{\rho^{(k)}} \right) \right]$\\

$=\frac{1}{h_r^\alpha} \displaystyle\sum_{k \in I_r} \left[ M_k \left( \frac{ \vert t_{km}(A_k(y)-L) \vert + \vert t_{km}(A_k(z)-L)\vert }{\rho^{(k)}} \right) \right]$\\

$\leq \frac{1}{h_r^\alpha} \displaystyle\sum_{k \in I_r} \left[ \frac{1}{2} M_k \left( \frac{ 2 \vert t_{km}(A_k(y)-L) \vert}{\rho^{(k)}} \right) + \frac{1}{2} M_k \left( \frac{ 2 \vert t_{km}(A_k(z)-L)\vert }{\rho^{(k)}} \right) \right]$\\

$< \frac{1}{2} \frac{1}{h_r^\alpha} \displaystyle\sum_{k \in I_r} K_1 \left(\frac{ \vert t_{km}(A_k(y)-L) \vert}{\rho^{(k)}} \right) M_k(2) + \frac{1}{2} \frac{1}{h_r^\alpha} \displaystyle\sum_{k \in I_r} K_2 \left(\frac{ \vert t_{km}(A_k(z)-L) \vert}{\rho^{(k)}} \right) M_k(2)$\\

$\leq \frac{1}{2} \frac{1}{h_r^\alpha} \displaystyle\sum_{k \in I_r} K_1 \left(\frac{ \vert t_{km}(A_k(y)-L) \vert}{\rho^{(k)}} \right) \mbox{sup} M_k(2) + \frac{1}{2} \frac{1}{h_r^\alpha} \displaystyle\sum_{k \in I_r} K_2 \left(\frac{ \vert t_{km}(A_k(z)-L) \vert}{\rho^{(k)}} \right) \mbox{sup} M_k(2)$\\

$\rightarrow 0 \mbox{~as~} r \rightarrow \infty.$\\

Hence $x \in \hat N_\theta^\alpha (A, \mathscr{M})$. This completes the proof.\\

{\bf{Theorem 3.6. }} Let $A=(a_{mk})$ be an infinite matrix of complex numbers and let $\mathscr{M}=(M_k)$ be a Musielak-Orlicz function satisfying $\Delta_2$-condition. If

$$\displaystyle\lim_{\nu \rightarrow \infty} \displaystyle\inf_k \frac{M_k\left(\frac{\nu}{\rho^{(k)}}\right)}{\frac{\nu}{\rho^{(k)}}} >0, \mbox{~for some ~} \rho^{(k)} >0,$$
then, $\hat N_\theta^\alpha(A)=\hat N_\theta^\alpha(A, \mathscr{M})$.

Proof: If  $\hat N_\theta^\alpha(A)=\hat N_\theta^\alpha(A, \mathscr{M})$ for some $\rho^{(k)}>0$, then there exists a number $\gamma>0$ such that

$$M_k\left(\frac{\nu}{\rho^{(k)}} \right) \geq \gamma \left( \frac{\nu}{\rho^{(k)}} \right), \forall  \nu >0, \mbox{~and some~} \rho^{(k)} > 0.$$

Let $x \in \hat N_\theta^\alpha(A, \mathscr{M})$ . Then,

\begin{eqnarray*}
\frac{1}{h_r^\alpha} \displaystyle\sum_{k \in I_r} \left[ M_k \left( \frac{\vert t_{km}(A_k(x)-L) \vert}{\rho^{(k)}} \right)\right]
& \geq & 
 \frac{1}{h_r^\alpha} \displaystyle\sum_{k \in I_r} \nu \left[ \gamma \left( \frac{t_{km}(A_k(x)-L) \vert}{\rho^{(k)}} \right)\right]\\ & = &
\gamma \frac{1}{h_r^\alpha}\displaystyle\sum_{k \in I_r} \left( \frac{\vert t_{km}(A_k(x)-L) \vert}{\rho^{(k)}} \right)
\end{eqnarray*}

Hence, $x \in \hat N_\theta^\alpha(A)$. This completes the proof.\\

{\bf{Theorem 3.7. }} Let $\mathscr{M}=(M_k)$ be a Musielak-Orlicz function where $(M_k)$ is pointwise convergent. Then, $\hat N_\theta^\alpha(A, \mathscr{M}, (s)) \subset \hat S_\theta^\alpha(A, \mathscr{M}, (s))$ if and only if $\displaystyle\lim_k M_k \left(\frac{ \nu}{\rho^{(k)}}\right) > 0$ for some $\nu>0, \rho^{(k)}>0$.\\

Proof : Let $\varepsilon >0$ and $x \in \hat N_\theta^\alpha(A, \mathscr{M}, (s))$. \\

Also, if $\displaystyle\lim_k M_k\left( \frac{\nu}{\rho^{(k)}}\right) >0$, then there exists a number $c > 0$ such that

$$M_k \left(\frac{\nu}{\rho^{(k)}}\right) \geq c, \mbox{~for~} \nu > \varepsilon.$$

Let us consider, $I_r^1=\left\{ i \in I_r: \left[ M_k \left(\frac{\vert t_{km}(A(x)-L) \vert}{\rho^{(k)}}\right)\right] \geq \varepsilon \right\}$. \\

Then, 
\begin{eqnarray*}
\frac{1}{h_r^\alpha} \displaystyle\sum_{k\in I_r} \left[ M_k \left(\frac{\vert t_{km}(A(x)-L) \vert}{\rho^{(k)}} \right)\right]^{s_k}
&\geq &
 \frac{1}{h_r^\alpha} \displaystyle\sum_{k\in I_r^1} \left[ M_k \left(\frac{\vert t_{km}(A(x)-L) \vert}{\rho^{(k)}} \right)\right]^{s_k}\\
&\geq &
c \frac{1}{h_r^\alpha} \vert t_{km}(A_0(\varepsilon) \vert
\end{eqnarray*}

Hence, it follows that $ x \in \hat S_\theta^\alpha(A, \mathscr{M}, (s))$.\\

Conversely, let us assume that the condition does not hold good. For a number $\nu>0$ , let $\displaystyle\lim_k M_k\left(\frac{\nu}{\rho^{(k)}}\right)=0$ for some $\rho>0$. Now, we select a lacunary sequence $\theta=(n_r)$ such that $M_k \left(\frac{\nu}{\rho^{(k)}}\right)< 2^{-r}$ for any $k > n_r$.\\

Let $A=I$ and define a sequence $x$ by putting,

\[ A_k(x) = \left\{ \begin{array}{ll}
\nu & \mbox{if $n_{r-1} < k \leq \frac{n_r+n_{r-1}}{2}$};\\
\theta & \mbox{if $\frac{n_r+n_{r-1}}{2} < k \leq n_r$}.\end{array} \right. \]

Therefore,

\begin{eqnarray*}
\frac{1}{h_r^\alpha}\displaystyle\sum_{k \in I_r}\left[ M_k \left( \frac{\vert A_k(x) \vert}{\rho^{(k)}} \right)\right]^{s_k}
&=&
\frac{1}{h_r^\alpha} \displaystyle\sum_{n_{r-1} < k \leq \frac{(n_r+n_{r-1})}{2}} M_k \left( \frac{\nu}{\rho^{(k)}} \right)\\
&<&
\frac{1}{h_r^\alpha} \frac{1}{2^{r-1}} \left[ \frac{n_r+n_{r-1}}{2} - n_{r-1} \right]\\
&=&
\frac{1}{2^r} \rightarrow 0 \mbox{~as~} r \rightarrow \infty.
\end{eqnarray*}

Thus we have $x \in \hat N_\theta^{\alpha 0} (A, \mathscr{M}, (s))$. \\

But, 
\begin{eqnarray*}
\displaystyle\lim_{r \rightarrow \infty} \frac{1}{h_r^\alpha} \left\vert \left\{ k\in  I_r: \displaystyle\sum_{k \in I_r}\left[M_k \left(\frac{\vert t_{km}(A(x))\vert}{\rho^{(k)}} \right)\right]^{s_k} \geq \varepsilon \right\} \right\vert
&=&
\displaystyle\lim_{r \rightarrow \infty} \frac{1}{h_r^\alpha} \left \vert \left \{ k \in \left( n_{r-1}, \frac{n_r + n_{r-1}}{2} \right) : \displaystyle\sum_{k \in I_r} \left[M_k \left( \frac{\nu}{\rho^{(k)}} \right)\right]^{s_k} \geq \varepsilon \right\} \right\vert\\
&=&
\displaystyle\lim_{r \rightarrow \infty} \frac{1}{h_r^\alpha} \frac{n_r-n_{r-1}}{2} \\
&=& 
\frac{1}{2}.
\end{eqnarray*}

So, $x \notin \hat S_\theta^\alpha(A, \mathscr{M}, (s))$.  \\

{\bf{Theorem 3.8. }} Let $\mathscr{M}=(M_k)$ be a Musielak-Orlicz function. Then $\hat S_\theta^\alpha(A, \mathscr{M},(s)) \subset \hat N_\theta^\alpha(A, \mathscr{M}, (s))$ if and only if $\displaystyle\sup_\nu \displaystyle\sup_k M_k \left(\frac{\nu}{\rho^{(k)}}\right) < \infty$.\\

Proof: Let $x \in \hat S_\theta^\alpha(A, \mathscr{M},(s))$. Suppose $h(\nu)=\displaystyle\sup_k M_k \left(\frac{\nu}{\rho^{(k)}} \right)$ and $h=\displaystyle\sup_\nu h(\nu)$. Let

$$I_r^2=\left\{k \in I_r: M_k \left(\frac{\vert t_{km}(A(x)-L) \vert}{\rho^{(k)}} \right) < \varepsilon \right\}.$$

Now, $M_k(\nu) \leq h$ for all $k, \nu >0$. So,
\begin{eqnarray*}
\frac{1}{h_r^\alpha} \displaystyle\sum_{k\in I_r} \left[ M_k \left(\frac{\vert t_{km}(A(x)-L) \vert}{\rho^{(k)}} \right)\right]^{s_k}
&= &
 \frac{1}{h_r^\alpha} \displaystyle\sum_{k\in I_r^1} \left[ M_k \left(\frac{\vert t_{km}(A(x)-L) \vert}{\rho^{(k)}} \right)\right]^{s_k}\\
&+&
 \frac{1}{h_r^\alpha} \displaystyle\sum_{k\in I_r^2} \left[ M_k \left(\frac{\vert t_{km}(A(x)-L) \vert}{\rho^{(k)}} \right)\right]^{s_k}\\
&\leq &
h \frac{1}{h_r^\alpha} \vert t_{km}(A_0(\varepsilon) \vert + h(\varepsilon).
\end{eqnarray*}

Hence, as $\varepsilon \rightarrow 0$, it follows that $x \in \hat N_\theta^\alpha(A, \mathscr{M}, (s))$.\\

Conversely, suppose that
$$\displaystyle\sup_\nu \displaystyle\sup_k M_k \left( \frac{\nu}{\rho^{(k)}} \right) =\infty.$$

Then, we have
$$0< \nu_1< \nu_2 <...< \nu_{r-1}< \nu_r<...$$

so that $M_{n_r} \left(\frac{\nu_r}{\rho^{(k)}} \right) \geq h_r^\alpha$ for $r \geq 1$. Let $A=I$. We set a sequence $x=(x_k)$ by,
\[ A_k(x) = \left\{ \begin{array}{ll}
\nu_r & \mbox{if $k=n_r$ for some $r=1,2,..$};\\
\theta & \mbox{otherwise}.\end{array} \right. \]

Then,

\begin{eqnarray*}
\displaystyle\lim_{r \rightarrow \infty} \frac{1}{h_r^\alpha} \left\vert \left\{ k \in I_r : \left[ \displaystyle\sum_{k \in I_r} M_k \left(\frac{\vert t_{km}(A_k(x)) \vert}{\rho^{(k)}} \right) \right] ^{s_k} \geq \varepsilon \right\} \right \vert 
&=&
\displaystyle\lim_{r \rightarrow \infty} \frac{1}{h_r^\alpha} \\
&=&
0
\end{eqnarray*}

Hence, $x \in \hat S_\theta^\alpha(A, \mathscr{M}, (s))$.\\

But,
\begin{eqnarray*}
\displaystyle\lim_{r \rightarrow \infty} \frac{1}{h_r^\alpha} \displaystyle\sum_{k \in I_r} \left[M_k\left( \frac{\vert t_{km} (A_k(x)-L) \vert}{\rho^{(k)}}\right)\right]
&=&
\displaystyle\lim_{r \rightarrow \infty} \frac{1}{h_r^\alpha} \left[M_{n_r}\left(\frac{\vert \nu_r- L\vert}{\rho^{(k)}}\right)\right]\\
&\geq &
\displaystyle\lim_{r \rightarrow \infty} \frac{1}{h_r^\alpha} h_r^\alpha\\
&=&1
\end{eqnarray*}

So, $x \in \hat N_\theta^\alpha(A, \mathscr{M}, (s))$.

\end{document}